\documentclass[12pt]{article}
\usepackage{amssymb,amsmath}
\def\OPLUS#1{\raisebox{-3pt}{\mbox{$\begin{array}{c}\mbox{\Large$\oplus$}\\[-6pt]\sc #1\end{array}$}}}
\def\hang{\hangindent\parindent}
\def\tex#1{\indent\llap{[#1]\enspace}\ignorespaces}
\def\re{\par\hang\tex}
\def\a{\alpha}
\def\b{\beta}

\def\SUM#1#2{\mbox{$\sum\limits_{#1}^{#2}$}}
\def\D{\Delta}
\def\UU{{\mathcal U}}

\def\g{\gamma}
\def\G{\Gamma}
\def\Der{{\rm Der}}
\def\Inn{{\rm Inn}}

\def\pp{{\mathcal P}}
\def\Ker{{\rm Ker}}
\def\es{\varepsilon}
\def\Im{{\rm Im}}
\def\span{{\rm span}}

\def\si{\sigma}
\def\v{\varphi}

\def\sc{\scriptstyle}
\def\ssc{\scriptscriptstyle}

\def\cl{\centerline}

\def\ol{\overline}

\def\wt{\widetilde}

\def\rar{\rightarrow}
\def\Rar{\Longrightarrow}

\def\PP{{\cal P}}
\def\bs{\backslash}

\def\vs{\vspace*}

\def\ni{\noindent}

\def\Z{\mathbb{Z}{\ssc\,}}

\def\F{\mathbb{F}{\ssc\,}}
\def\QED{\hfill$\Box$}

\numberwithin{equation}{section}
\newtheorem{theo}{Theorem}[section]
\newtheorem{defi}[theo]{Definition}

\newtheorem{lemm}[theo]{Lemma}
\newtheorem{clai}{Claim}
\newtheorem{subc}{Subclaim}
\def\adddot{ \ }
\def\addbra{$\!\!\!${\bf)}\ \ }
\def\PP{{\mathcal H}}
\def\pp{\PP}
\begin{document}
\cl{{\Large \bf Hamiltonian type Lie bialgebras}\footnote {Supported
by NSF grant 10471096 of China and ``One Hundred Talents Program''
from University of Science and Technology of China
} } \vs{6pt}

\cl{Bin Xin$^{\,1)}$,  Guang'ai Song$^{\,2)}$, Yucai Su$^{\,3)}$}

 \cl{\small
$^{1)}$Department of Mathematics, Shanghai Jiao Tong University,
 Shanghai 200240, China}

\cl{\small $^{2)}$College of Mathematics and Information
\vs{-3pt}Science, Shandong Institute of Business and} \cl{\small
Technology, Yantai, Shandong 264005, China}

 \cl{\small \small $^{3)}$Department of
Mathematics, University of Science and Technology of \vs{-3pt}China}
\cl{\small Hefei 230026, China}

 \cl{\small E-mail:  ycsu@ustc.edu.cn} \vs{6pt}

{\small
\parskip .005 truein
\baselineskip 3pt \lineskip 3pt

\noindent{{\bf Abstract.} We first prove that, for any generalized
Hamiltonian type Lie algebra $\PP$, the first cohomology group
$H^1(\PP,\PP \otimes\PP)$ is trivial. We then show that all Lie
bialgebra structures on  $\PP$ are triangular. \vs{5pt}

}

}
\parskip .001 truein\baselineskip 6pt \lineskip 6pt

\vs{6pt}

\cl{\bf\S1. \
Introduction}\setcounter{section}{1}\setcounter{equation}{0} To
study quantum groups, Drinfeld [D1] (cf.~[D2]) introduced the notion
of Lie bialgebras. We start by recalling the definition. Let $L$ be
a vector space over a field $\F$ of characteristic zero. Denote by
$\tau$ the {\it twist map} of $L \otimes L$, namely, $\tau(x\otimes
y)= y \otimes x$ for $x,y\in L.$ Denote by $\xi:L\otimes L\otimes
L\to L\otimes L\otimes L$ the linear map which cyclicly permutes the
coordinates, i.e.,
$$ \xi=(1\otimes\tau)(\tau\otimes 1):x_{1} \otimes x_{2} \otimes
x_{3}\mapsto x_{2} \otimes x_{3} \otimes x_{1}\mbox{ \ for \
}x_1,x_2,x_3\in L,$$ where (more generally, throughout this paper)
$1$ denotes the {\it identity map}. Then a {\it Lie algebra} is a
pair $(L,\v)$ of a vector space $L$ and a linear map $\v :L \otimes
L \rar L$ (the {\it bracket} of $L$) satisfying the following
conditions:
\begin{eqnarray*}
\!\!\!\!\!\!\!\!\!\!\!\!&&
\Ker(1\otimes1-\tau) \subset \Ker\,\v \mbox{ \ ({\it strong skew-symmetry}), \ \ and}\\
\!\!\!\!\!\!\!\!\!\!\!\!&& \v \cdot (1 \otimes \v ) \cdot
(1\otimes1\otimes1 + \xi +\xi^{2}) =0 : \ L \otimes L \otimes L \rar
L\mbox{ \ ({\it Jacobi identity}).}
\end{eqnarray*}
Dually, a {\it Lie coalgebra} is a pair $(L, \D)$ of a vector space
$L$ and a linear map $\D: L \to L \otimes L$ ({\it cobracket} of
$L$) satisfying the following conditions:
\begin{eqnarray}
\label{cLie-s-s} \!\!\!\!\!\!\!\!\!\!\!\!&&
\Im\,\D \subset \Im(1\otimes1- \tau) \mbox{ \ ({\it anti-commutativity}), \ \ and}\\
\label{cLie-j-i} \!\!\!\!\!\!\!\!\!\!\!\!&& (1\otimes1\otimes1 + \xi
+\xi^{2}) \cdot (1 \otimes \D) \cdot \D =0:\ L \to L \otimes L
\otimes L\mbox{ \ ({\it Jacobi identity}).}
\end{eqnarray}
\begin{defi} \
\rm  A {\it Lie bialgebra} is a triple $(L,\v, \D )$, where $(L,
\v)$ is a Lie algebra and $(L, \D)$ is a Lie coalgebra satisfying
the following compatibility condition:
\begin{eqnarray}
\label{bLie-d} \!\!\!\!\!\!\!\!\!\!\!\!&& \mbox{$\D\cdot  \v (x, y)
= x \cdot \D y - y \cdot \D x$ \ for \ $x, y \in L$,}
\end{eqnarray}
where the symbol ``$\cdot$'' on the right-hand side denotes the {\it
adjoint diagonal action} of $L$ on $L\otimes L$ under which, for all
$x,a_i,b_i\in L$,
\begin{equation}
\label{diag} \mbox{$x\cdot (\sum\limits_{i} {a_{i} \otimes b_{i}}) =
\sum\limits_{i} ( {[x, a_{i}] \otimes b_{i} + a_{i} \otimes [x,
b_{i}]})$,}
\end{equation}
where, in general $[x,y]=\v(x\otimes y)$ for $x,y \in L$.
\end{defi}
%

\begin{defi}
\label{def2} \rm \vskip-6pt\begin{itemize}\parskip-5pt
\item[{(1)}]
 A {\it coboundary Lie bialgebra} is a $4$-tuple $(L, \v, \D,r),$ where
$(L, \v, \D)$ is a Lie bialgebra, $r \in \Im(1\otimes1 - \tau)
\subset L \otimes L$, and $\D=\D_r$ where $\D_r$ (the {\it
coboundary} of $r$) is defined by
\begin{equation}
\label{D-r}
\D_r (x) = x \cdot r \mbox{ \ \ for \ \ }x \in L.
\end{equation}
\item[(2)] A coboundary Lie bialgebra $(L, \v,\D, r)$
is called {\it triangular} if it satisfies the {\it classical
triangle equation} also known as the {\it classical Yang-Baxter
Equation} (CYBE)
\begin{equation}
\label{CYBE}
c(r):= [r^{12} , r^{13}] +[r^{12} , r^{23}] +[r^{13} , r^{23}]=0,
\end{equation}
 where the $r^{ij}$ are defined as follows: If $r =\sum_{i} {a_{i} \otimes
b_{i}} \in L \otimes L$, then (upon identifying an element of $L$
with its image in $\UU(L)$, the {\it universal enveloping algebra}
of $L$)
$$\begin{array}{l}
 r^{12} = \sum \limits_{i}{a_{i} \otimes b_{i}
\otimes 1}=r\otimes1, \\[12pt]
r^{13}= \sum \limits_{i} {a_{i} \otimes 1 \otimes b_{i}}=(1\otimes\tau)(r\otimes1)=(\tau\otimes1)(1\otimes r), \\[12pt]
r^{23} = \sum \limits_{i}{1 \otimes a_{i} \otimes b_{i}}=1\otimes r,
\end{array}$$
are elements of ${\UU (L) \otimes \UU (L) \otimes \UU (L)}$, $1$ is
the identity element of $\UU (L)$ and the above-defined brackets are
just the {\it commutators} of the given elements in the associative
unitary algebra ${\UU (L) \otimes \UU (L) \otimes \UU (L)}$.
\end{itemize}
\end{defi}

During the past two decades, many papers on Lie bialgebras have
appeared among them for examples, those listed in our references,
e.g., [D1, D2, Dzh, M, NT, N, S, SS, T, WSS]. Lie bialgebras arise
naturally in the study of Hamiltonian mechanics and Poisson Lie
group. Michaelis [M] presented a class of infinite-dimensional Lie
bialgebras containing the Virasoro and Witt algebras. Such Lie
bialgebras were subsequently classified by Ng and Taft [NT] (cf.~[N,
T]). Song and Su [SS] (see also [S]), considered Lie bialgebra
structures on the graded Lie algebras of generalized Witt type
constructed in [DZ]. In the present paper, we consider Lie bialgebra
structures on Lie algebras of generalized  Hamiltonian type (or
Cartan type $H$) defined in [OZ] (see also [X]).
Such Lie algebras attracted our attention because they have a
Poisson algebra structure (cf.~(\ref{PS--})). Poisson algebras are
fundamental algebraic structures on the phase spaces in classical
mechanics; they are also of major interest in symplectic geometry.

In what follows, let $\G$ be any {\it nondegenerate additive
subgroup} of $\F^{2n}$ (by which we mean that  $\G$ contains an
$\F$-basis of $\F^{2n}$). Choose an $\F$-basis $\{\es_p\,|\,1\le
p\le 2n\}\subset\G$ of $\F^{2n}$. Any element $\a\in\G$ can be
written as
\begin{eqnarray}
\label{alpha} \!\!\!\!\!\!\!\!&&\a=(\a_1,\a_{\ol1},...,\a_n,\a_{\ol
n})\in\G\mbox{ \ with \ }\a_1,\a_{\ol1},...,\a_n,\a_{\ol n}\in\F
\mbox{ and }\ol p=n+p.
\end{eqnarray}
Set $$\si_p=\si_{\ol p}=\es_p+\es_{\ol p}~\mbox{ \ for \ }1\le p\le
n.$$ Let $\ol{\pp}=\ol{\pp}(2n,\G):={\rm
span}{\{t^{\a}\,|\,\a\in\G\}}$ be the {\it group algebra} with
product given by $ t^{\a}\cdot t^{\b}= t^{\a+\b}\:$. Define
$[\cdot,\cdot]$ as follows,
\begin{eqnarray}
\label{braket}
\!\!\!\!\!\!\!\!&&[t^{\a},t^{\b}]=\SUM{i=1}{n}(\a_i\b_{\ol
i}-\b_i\a_{\ol i})t^{\a+\b-\si_i}\mbox{ \ for \ } \a,\b\in \G.
\end{eqnarray}
Then it can be shown (cf.~[X, SX]) that
$({\ol\pp,\cdot,[\cdot,\cdot]})$ is a {\it Poisson algebra\,}
satisfying the following compatible condition:
\begin{equation}\label{PS--} [u,v\cdot w]=[u,v]\cdot w+v\cdot[u,w]\mbox{ \
\ for \ }u,v,w\in\ol\pp.\end{equation} Let $\pp=\ol{\pp}/\F\cdot1$
where $1=t^0$. Then $\pp$ is a simple Lie algebra called a {\it
Hamiltonian Lie algebra of Cartan type} [OZ, X] (see also [Su, SX]).

Let $V$ be a module over $\pp$. Denote by $\Der(\pp,V)$ the set of
\textit{derivations} $d:\pp\to V$ satisfying
\begin{equation}
\label{deriv} d([x,y])=x\cdot d(y)-y\cdot d(x)\mbox{ \ \ for \
}x,y\in \pp,
\end{equation}
and denote by $\Inn(\pp,V)$ the set consisting of the \textit{inner
derivations} $a_{\rm inn}$ defined, for each $a\in V$, by
\begin{equation}
\label{inn} a_{\rm inn}:x\mapsto x\cdot a\mbox{ \ \ for \ }x\in \pp.
\end{equation}
It is well-known that the first cohomology group $H^1(\pp,V)$ of
$\pp$ with coefficients in $V$ is isomorphic to
$\Der(\pp,V)/\Inn(\pp,V)$, i.e.,
$H^1(\pp,V)\cong\Der(\pp,V)/\Inn(\pp,V).$

The main results of this paper is the following.
\begin{theo}
\label{main} \vskip-6pt\begin{itemize}\parskip-5pt \item[{\rm(1)}]
Every Lie bialgebra structure on $\pp$ is  triangular.
\item[{\rm(2)}] An element $r\in\pp\otimes\pp$ satisfies CYBE $(\ref{CYBE})$
if and only if $c(r)$ is \textbf{ad-invariant}, i.e.,
\begin{equation}
\label{MYBE} x\cdot c(r)=0\mbox{ \ for all \ }x\in \pp.
\end{equation}
\item[{\rm(3)}] Regarding $V=\pp\otimes\pp$ as an $\pp$-module under the
adjoint diagonal action $(\ref{diag})$ of $\pp$, one has that
$H^1(\pp,V)=\Der(\pp,V)/\Inn(\pp,V)=0$. \end{itemize}
\end{theo}
%

 \vskip10pt

\cl{\bf\S2. \ Proof of the main results}\setcounter{section}{2}
\setcounter{theo}{0}\setcounter{equation}{0}

First we retrieve some useful results from Drinfeld [D1, D2],
Michaelis [M], Ng-Taft [NT] and combine them into the following
theorem.

\begin{theo}
\label{some}
 \vskip-6pt\begin{itemize}\parskip-5pt \item[{\rm(1)}]
For a Lie algebra $L$ and $r\in \Im(1\otimes1-\tau)\subset L\otimes
L,\,\D=\D_r$, one has that
\begin{equation}
\label{add-c} (1 + \xi + \xi^{2}) \cdot (1 \otimes \D) \cdot \D (x)
= x \cdot c (r) \mbox{ \ for all \ }x\in L.
\end{equation}
Thus  the triple $(L,[\cdot,\cdot], \D_r)$ is a Lie bialgebra if and
only if $r$ satisfies $(\ref{MYBE})$.
\item[{\rm(2)}]
Let $L$ be a Lie algebra containing two elements $a,b$ satisfying
$[a, b] =b$, and set $r=a \otimes b -b \otimes a$. Then $\D_{r}$
equips $L$ with the structure of a triangular Lie bialgebra.
\end{itemize}
\end{theo}

We shall divide the proof of Theorem \ref{main} into several lemmas.
First we need some preparations. Define the linear map
$\pi:\G\to\F^{n}$ by
\begin{equation}
\label{ptl} \pi(\a)=(\a_1-\a_{\ol1},\a_2-\a_{\ol2},...,\a_n-\a_{\ol
n})\in \F^{n}\mbox{ \ \ for \ }\a\in\G.
\end{equation}
Let $$G:=\pi(\G)=\{\pi(\a)\,|\,\a\in\G\}\subset\F^n.$$ We always
denote an element $\mu$ of $\F^n$ as $$\mu=(\mu_1,...,\mu_n).$$ Then
$\pp$ is a $G$-graded Lie algebra (but in general not
finitely-graded):
\begin{equation}
\label{pp} \pp=\OPLUS{\sc \mu\in G}\pp_{\mu},\mbox{ \ \ where \
}\pp_{\mu} ={\rm span}\{t^{\b}\,|\,\pi(\b)=\mu\}.
\end{equation}
Also $$V:=\pp\otimes\pp=\OPLUS{\mu\in G}V_\mu\mbox{ \ is $G$-graded,
with } V_\mu=\mbox{$\SUM{\nu+\lambda=\mu}{}$}\pp_\nu\otimes
\pp_\lambda.$$
 It is easily seen from (\ref{braket}) that
\begin{equation}\label{e-t-si} t^{\si_p}|_{V_\mu}=\mu_p\cdot 1_{V_\mu}
\mbox{ \ \ for \ }\mu\in G,\ \ p=1,2,...,n.
\end{equation}
For $1\le p\le n$, we denote
\begin{equation}\label{e-p-p}
\pp^p=\span\{t^\a\,|\,\a_p=\a_{\ol p}=0\},\ \ \ \pp^{\le
p}=\span\{t^\a\,|\,\a_q=\a_{\ol q}=0,\,q\le p\}. \end{equation}
Obviously, $\pp^{\le n}=0$. Set $$V^{p}=\pp^p\otimes
\pp+\pp\otimes\pp^p,\,\ \ V^{\le p}=\pp^{\le p}\otimes
\pp+\pp\otimes \pp^{\le p}.$$
\begin{lemm}\adddot
\label{lemma2.3} Let $v\in V,\,1\le p\le n$. Suppose
$t^{k\es_p}\cdot v=t^{k\es_{\ol p}}\cdot v=0$ for $k\in\Z$. Then
$v\in V^{p}$.
\end{lemm}
\ni{\it Proof.} Write $v=\sum_{\a,\b\in\G}c_{\a,\b}t^{\a}\otimes
t^{\b}$ for some $c_{\a,\b}\in\F$, where
$\{(\a,\b)\,|\,c_{\a,\b}\ne0\}$ is a finite set. Choose $k\gg0$ such
that $c_{\a+k\es_p-\si_p,\b-k\es_p+\si_p}=0$ for all $\a,\b$. Then
$$
\begin{array}{ll}
0\!\!\!&=t^{k\es_p}\cdot v
=k\sum\limits_{\a,\b\in\G}c_{\a,\b}(\a_{\ol p}
t^{\a+k\es_p-\si_p}\otimes t^\b +\b_{\ol b}t^\a\otimes
t^{\b+k\es_p-\si_p})
\\[12pt]
&=k\sum\limits_{\a,\b\in\G}(\a_{\ol p} c_{\a,\b}+(\b_{\ol
p}+1)c_{\a+k\es_p-\si_p,\b-k\es_p+\si_p}) t^{\a+k\es_p-\si_p}\otimes
t^\b
\\[12pt]
&= k\sum\limits_{\a,\b\in\G}\a_{\ol p}c_{\a,\b}
t^{\a+k\es_p-\si_p}\otimes t^\b. \end{array}$$ Thus
$c_{\a,\b}\ne0\Rar\a_p=0$. Similarly $c_{\a,\b}\ne0\Rar\b_p=0$.
Using $\ol p$ instead of $p$, we have $c_{\a,\b}\ne0\Rar\a_{\ol
p}=\b_{\ol p}=0$. \QED\vskip5pt

\begin{lemm}\adddot
\label{last} Suppose $r\in V$ such that $a\cdot r\in {\rm
Im}(1\otimes1-\tau)$ for all $a\in\PP$. Then $r\in{\rm
Im}(1\otimes1-\tau)$.
\end{lemm}
\ni{\it Proof.} (cf.~[SS]) Write $r=\sum_{\a\in\G}r_\a$ with
$r_\a\in V_\a$. Obviously, $r\in{\rm Im}(1\otimes1-\tau)$ if and
only if $r_\a\in{\rm Im}(1\otimes1-\tau)$ for all $\a\in\G.$ Thus
without loss of generality, we can suppose $r=r_\a$ is
homogeneous. Note that $\PP\cdot {\rm
Im}(1\otimes1-\tau)\subset{\rm Im}(1\otimes1-\tau)$. We shall
subtract $r_\a$ by $r_\a-u$ for some $u\in{\rm
Im}(1\otimes1-\tau)$ to reduce it into $0$.

First assume $\a\ne0$. By (\ref{e-t-si}), we can choose $t^{\si_p}$
such that $t^{\si_p}|_{V_\a}=\a_p\cdot 1_{V_\a}\ne0$. Thus,
$r_\a=\a_p^{-1}t^{\si_p}\cdot r_\a\in {\rm Im}(1\otimes1-\tau)$.

Assume now $\a=0$. Write
$r_0=\Sigma_{\g,\b\in\G}\,c_{\g,\b}t^{-\b}\otimes t^{\g+\b}$, where
$\pi(\g)=0$. Choose a total order on $\G$ compatible with its group
structure. Since $u_{\g,\b}:=t^{-\b}\otimes
t^{\g+\b}-t^{\g+\b}\otimes t^{-\b}\in {\rm Im}(1\otimes1-\tau)$, by
replacing $r$ by $r-u$, where $u$ is some combination of
$u_{\g,\b}$, we can suppose
\begin{equation}
\label{last-c} c_{\g,\b}\ne0\ \Longrightarrow\ -\b\preceq \g+\b\,.
\end{equation}
Now assume that $c_{\g,\b}\ne0$ for some $-\b\preceq\g+\b$. Suppose
$\b_{\ol i}\neq0$ for some $i$. Taking $m\gg0$, we get
$$
\begin{array}{ll}
\!\!\!&t^{m\es_i}\cdot r_0
=-m\sum\limits_{\g,\b\in\G}c_{\g,\b}(\b_{\ol
i}t^{-\b+m\es_i-\si_i}\otimes t^{\g+\b} -(\g_{\ol i}+\b_{\ol
i})t^{-\b}\otimes t^{\g+\b+m\es_i-\si_i})\,.
\\[12pt]
\end{array}$$
From this, we obtain $t^{m\es_i}\cdot r_0\notin {\rm
Im}(1\otimes1-\tau).$ Thus $\b_{\ol i}=0$ for $1\le i\le n$. By
symmetry, $\b_i=0$ for $1\le i\le n$. Similarly, $\g_i=0$ for $1\le
i\le 2n$. Namely, $c_{\g,\b}\ne0$ implies $\g=\b=0$. This proves
$r_0=0$. \QED \vskip6pt
 Theorem \ref{main}(2) follows from the following more general
result.

\begin{lemm}\adddot
\label{lemma2} Regarding $\pp[m]=\pp\otimes\cdots\otimes\pp$
$({\ssc\,}m$ copies$)$ as a $\pp$-module under diagonal action,
suppose $c\in\pp[m]$ satisfying $a\cdot c=0$ for all $a\in\pp$. Then
$c=0$.
\end{lemm}
\ni{\it Proof.} Similar to the proof of Lemma \ref{lemma2.3}, we can
obtain $c\in \pp[m]^{\le n}=\{0\}$ (cf.~(\ref{e-p-p})). \QED
\vskip7pt
%
Theorem \ref{main}(3) follows from the following most technical
theorem in this paper.
\begin{theo}\adddot
\label{lemma3?} $\Der(\pp,V)=\Inn(\pp,V)$, where $V=\pp\otimes\pp$.
\end{theo}
\ni{\it Proof.} First we would like to remark that although we have
observed (cf.~[SS, WSS]) that the analogous result of the theorem
holds for many Lie algebras, we are unfortunately unable to find a
unified way to prove this result for a general Lie algebra.

 We shall prove the theorem by several claims.
A derivation $d\in\Der(\pp,V)$ is {\it homogeneous of degree $\mu\in
G$} if $d(\pp_\nu) \subset V_{\mu +\nu}$ for $\nu \in G$. Denote
$$\Der(\pp, V)_\mu = \{d\in \Der(\pp, V) \,|\,{\rm deg\,}d =\mu\}.$$
\begin{clai}\adddot
\rm Let $d\in\Der(\pp,V)$. Then
\begin{equation}
\label{summable} \mbox{$d = \sum\limits_{\mu \in G} d_\mu , \mbox{ \
where \ }d_\mu \in \Der(\pp, V)_\mu,$}
\end{equation}
which holds in the sense that for every $u \in\pp$, only finitely
many $d_\mu (u)\neq 0,$ and $d(u) = \sum_{\mu \in G} d_\mu(u)$ (we
call such a sum in (\ref{summable}) {\it summable}).
\end{clai}

For $\mu\in G$, we define $d_\mu$ as follows:  For $u\in\pp_\nu$ and
$\nu\in G$, write $d(u)=\sum_{\g\in G}v_\g\in V$ with $v_\g\in
V_\g$, then we set $d_\mu(u)=v_{\mu+\nu}$. Obviously
$d_\mu\in\Der(\pp,V)_\mu$ and we have (\ref{summable}).
\begin{clai}\adddot
\rm If $0\ne \mu\in G$, then $d_\mu\in\Inn(\pp,V)$.
\end{clai}

Since $\mu\neq0$, we can choose $p$ such that $\mu_p\neq0$. For
$x\in\pp_\eta,\,\eta\in G$, applying $d_\mu$ to
$[t^{\si_p},x]=\eta_px$, by (\ref{deriv}), we have
\begin{equation*}
-x\cdot d_{\mu}(t^{\si_p})+t^{\si_p}\cdot
d_{\mu}(x)=\eta_pd_{\mu}(x).
\end{equation*}
Since $d_{\mu}(x)\in V_{\mu+\eta}$, we have $t^{\si_p}\cdot
d_{\mu}(x)=(\mu_p+\eta_p)d_{\mu}(x)$. Thus
$d_{\mu}(x)=u_p^{-1}x\cdot d_{\mu}(t^{\si_p})$. Namely,
$d_\mu=a_{\rm inn}$ for $a=u^{-1}_pd_\mu(t^{\si_p})$.

\begin{clai} \
\label{clai-sub1} \rm $d_0(t^{\si_p})=0$ for $1\le p\le n$.
\end{clai}

Applying $d_0$ to $[t^{\si_p},t^\a]=\pi(\a)_pt^\a$ for $\a\in\G,$ we
get
\begin{equation}
\label{t-p-0} -t^\a\cdot d_0(t^{\si_p})+\pi(\a)_{p\ssc\,}
t^\a=-t^\a\cdot d_0(t^{\si_p})+t^{\si_p}\cdot
d_0(t^\a)=\pi(\a)_{p\ssc\,}t^\a.
\end{equation}
That is, $t^\a\cdot d_0(t^{\si_p})=0.$ By Lemma \ref{lemma2},
$d_0(t^{\si_p})=0$.

\begin{clai}\adddot
\label{claim3} \rm $d_0\in\Inn(\pp,V)$.
\end{clai}

We shall prove by induction on $p=0,1,...,n,$ that by replacing
$d_0$ by $d_0-u_{\rm inn}$ for some $u\in V_0$,
\begin{equation}\label{e-ind}
d_0(x)\equiv0\ ({\rm mod\,}V^{\le p})\mbox{ \ \ for \ }x\in\pp.
\end{equation}
There is nothing to prove if $p=0$. Assume that $i>0$ and that
(\ref{e-ind}) holds for $p=i-1$. Now suppose $p=i$.
 The proof of
(\ref{e-ind}) will be done by several subclaims.

\begin{subc}\addbra
\label{sub2} \rm We can suppose $d_0(t^{\es_i})\equiv
d_0(t^{\es_{\ol i}})\equiv0\,({\rm mod\,}V^{\le i})$.
\end{subc}\ \indent

For any element $\g\in\G$, we can uniquely write (recall that $\pp$
is $G$-graded not $\G$-graded)
\begin{equation}
\label{s-form}
d_0(t^\g)=\SUM{\b,\a\in\G:\,\pi(\a)=0}{}c^\g_{\b,\a}t^{-\b}\otimes
t^{\b+\a+\g}\mbox{ \ (a finite sum),}
\end{equation}
for some $c^\g_{\b,\a}\in \F^*=\F\bs\{0\}$. For $\a\in\pi^{-1}(0)$,
define the linear map $D^\a$ by
$$D^\a(t^\g)=\SUM{\b\in\G}{}\,c_{\b,\a}^\g t^{-\b}\otimes t^{\b+\a+\g}.$$
It is easy to verify that $D^\a\in\Der(\pp,V)$ and
$d_0=\sum_{\a\in\pi^{-1}(0)}D^\a$.

 We shall prove Subclaim \ref{sub2})
in two steps.

\vskip5pt {\it Step 1.}~~We want to prove $d_0(t^{\es_i})\equiv
0\,({\rm mod\,}V^{\le i})$.\\

For $m,{m'}\in\Z$, we can write
\begin{equation}\begin{array}{l}
\label{d-t-i} D^{\a}(t^{m\es_i+{m'}\es_{\ol
i}})=\sum\limits_{\b\in\ol\G,\,j,k\in\Z}c_{\b,j,k}^{m,{m'}}
t^{-(\b+j\es_i+k\es_{\ol i})}\otimes t^{\b+j\es_i+k\es_{\ol
 i}+\a+m\es_i+{m'}\es_{\ol i}}\mbox{~\ (a finite sum),}
\end{array}\end{equation} where $\ol\G=\G\diagup(\Z\es_i+\Z\es_{\ol i})$,
and we can require that \begin{equation}\label{e-require} \b_p\in\Z\
\ \ \Longrightarrow\ \ \ \b_p=0\mbox{ \  for \ }p=i,\,\ol
i.\end{equation} First we consider $D^\a(t^{\es_i})$. Note that
\begin{eqnarray}
\label{-inn}
\!\!\!\!\!\!\!\!\!\!\!\!\!\!\!\!\!\![t^{\es_i},t^{-(\b+j\es_i+k\es_{\ol
i})}\otimes t^{\b+j\es_i+k\es_{\ol i}+\a+\si_i}]
&\!\!\!\!\!=\!\!\!\!\!&-(\b_{\ol i}+k)t^{-(\b+j\es_i+(k+1)\es_{\ol
i})}\otimes
t^{\b+j\es_i+(k+1)\es_{\ol i}+\a+\es_i} \nonumber\\
&&+(\b_{\ol i}+k+\a_{\ol i}+1)t^{-(\b+j\es_i+k\es_{\ol i})}\otimes
t^{\b+j\es_i+k\es_{\ol i}+\a+\es_i}.
\end{eqnarray}
From (\ref{-inn}), one can see that by replacing $D^\a$ by
$D^\a-u_{inn},$ where $u$ is a combination of some
$t^{-(\b+j\es_i+k\es_{\ol i})}\otimes t^{\b+j\es_i+k\es_{\ol
i}+\a+\si_i}$, and re-denoting $b_{\b,j,k}=c^{1,0}_{\b,j,k}$ (in
order to avoid confusion), we can suppose
\begin{equation}
\label{e-case-}
 ~b_{\b,j,k}\neq0\ \ \Longrightarrow\ \
 \left\{
\begin{array}{lll}
\mbox{1)}~~\b_{\ol i}\neq0,~ \b_{\ol i}+\a_{\ol i}\notin \Z,
&\!\!\mbox{and~~}k=0,\ \ \mbox{or}
\\[6pt]
\mbox{2)}~~\b_{\ol i}=0,~\a_{\ol i}\notin \Z,&\!\!\mbox{and}~~k=1 ,\
\ \mbox{or}
\\[6pt]
3)~~\b_{\ol i}\neq0,~ \b_{\ol i}+\a_{\ol i}\in \Z,&\!\!\mbox{and}
~~k=-(\b_{\ol i}+\a_{\ol i}+1)
,\ \ \mbox{or}\\[6pt]
4)~~\b_{\ol i}=0,~\a_{\ol i}\in \Z,&\!\!\mbox{and}~~k=1
~\mbox{or}~k=-(\a_{\ol i}+1).
\end{array}
\right. \end{equation} Applying $D^\a$ to $[t^{\es_i},t^{\es_{\ol
i}}]=0$, and by definition (\ref{d-t-i}), we obtain (here, the
summand is over the set as in (\ref{d-t-i}))
\begin{eqnarray}
 \label{d-0-i-i}
\!\!\!\!\!\!\!\!\!\!\!\!\!\!\!\!\!\!&&
\SUM{\b,j,k}{}(-(\b_i+j-1)b_{\b,j-1,k}+(\b_i+\a_i+j+1)b_{\b,j,k})t^{-(\b+j\es_i+k\es_{\ol
i})}\otimes t^{\b+j\es_i+k\es_{\ol i}+\a}
\nonumber\\
\!\!\!\!\!\!\!\!\!\!\!\!\!\!\!\!\!\!&& =\SUM{\b,j,k}{}((\b_{\ol
i}+k-1)c^{0,1}_{\b,j,k-1}-(\b_{\ol i}+\a_{\ol
i}+k+1)c^{0,1}_{\b,j,k}) t^{-(\b+j\es_i+k\es_{\ol i})}\otimes
t^{\b+j\es_i+k\es_{\ol i}+\a}.
\end{eqnarray}
 Fixing $j$, we obtain
\begin{equation}
\label{i-i}
(\b_i+\a_i+j+1)b_{\b,j,k}-(\b_i+j-1)b_{\b,j-1,k}=-(\b_{\ol
i}+\a_{\ol i}+k+1)c^{0,1}_{\b,j,k}+(\b_{\ol
i}+k-1)c^{0,1}_{\b,j,k-1}.
\end{equation}
Define $k_1\geq k_2$ to be the largest and the smallest integers
(depending on $j$) such that $c^{0,1}_{\b,j,k}\neq0,$ i.e,
\begin{equation}\label{e-define-k}k_1\geq k_2
\mbox{ \ \ and  \ \ }\left\{
\begin{array}{ll}
c^{0,1}_{\b,j,k}=0&\mbox{if \ }k>k_1\mbox{ or }k<k_2,\\[6pt]
c^{0,1}_{\b,j,k}\neq0&\mbox{if \ }k=k_1\mbox{ or }k=k_2.
\end{array}\right.
\end{equation}
Now assume that $b_{\b,j,k}\ne0$. We discuss (\ref{e-case-}) case by
case. \vskip5pt {\it Case 1})~~{$\b_{\ol i}\neq0,\, \b_{\ol
i}+\a_{\ol i}\notin \Z$ and $k=0$.}~~First suppose $k_1\ne-1$. In
(\ref{i-i}), replacing $k$ by $k_1+1$, the left-hand side becomes
zero, and we get $\b_{\ol i}+k_1=0$, which means that $\b_{\ol
i}\in\Z$. This and the assumption that $\b_{\ol i}\neq0$ contradict
(\ref{e-require}).  Thus $k_1=-1$.

Next suppose $k_2\ne0$. In (\ref{i-i}), replacing $k$ by $k_2$, we
get $\b_{\ol i}+k_2+\a_{\ol i}+1=0$, contradicting the assumption
that $\b_{\ol i}+\a_{\ol i}\notin \Z$. Thus $k_2=0$. But then
$k_1<k_2$, again a contradiction. Hence this case does not occur.
\vskip5pt {\it Case 2})~~{$\b_{\ol i}=0,\,\b_{\ol i}+\a_{\ol
i}\notin \Z$ and $k=1$.}~~Similar to Case 1), in (\ref{i-i}), by
replacing $k$ by $k_1+1$ and $k_2$ respectively, we obtain a
contradiction that $k_1=-\b_{\ol i}< k_2=1-\b_{\ol i}$. So this case
does not occur either.

\vskip5pt {\it Case 3})~~{$\b_{\ol i}\neq0,\,\a_{\ol i}\in \Z$ and
$k=-(\b_{\ol i}+\a_{\ol i}+1)$.}~~Using the same way, we obtain the
contradiction that $k_1=-(\b_{\ol i}+\a_{\ol i}+2)< k_2=-(\b_{\ol
i}+\a_{\ol i}+1)$.

The above 3 cases shows that for any $j$,
\begin{equation}\label{e-b-cond}
b_{\b,j,k}\ne0\ \ \ \Longrightarrow\ \ \ \b_{\ol i}=0,\,\a_i\in
\Z\mbox{ \ and \ $k=1$, \ or \ $k=-(\a_{\ol i}+1)$}. \end{equation}
Thus it remains to consider the following case.

\vskip5pt {\it Case 4})~~{$\b_{\ol i}=0,\,\a_{\ol i}\in \Z$ and
$k=1,$ or $k=-(\a_{\ol i}+1)$.}~~This case is more complicated, so
we shall discuss it in four subcases. \vskip5pt {\it Subcase
4.1})~~{$~\a_{\ol i}=-2$}. In this subcase, $b_{\b,j,k}\ne0$ implies
$k=1$. Then when $k$ takes value $1$ in (\ref{i-i}), the right-hand
side vanishes. We obtain (note that $\a_i=\a_{\ol i}=-2$)
\begin{equation}
\label{i-i-0} -(\b_i+j-1)b_{\b,j,1}=(\b_i+j-1)b_{\b,j-1,1}.
\end{equation}
If $b_{\b,j,1}\ne0$ for some $\b_i\notin\Z$, then (\ref{i-i-0})
implies $b_{\b,j',1}=b_{\b,j,1}\ne0$ for all $j'\in\Z$, which
contradicts the fact that $D^{\a}(t^{\es_i})$ is a finite sum. So
$\b_i\in\Z$, i.e., $\b_i=0$ by (\ref{e-require}).

If $b_{\b,j,1}\ne0$ for some $j>1$, then (\ref{i-i-0}) gives
$b_{\b,j+n,1}=b_{\b,j,1}\ne0$ for all $n\in \Z_+$, again a
contradiction. Similarly, if $b_{\b,j,1}\ne0$ for some $j<1$, we
have a contradiction. Thus $j=1$.

So suppose $b_{\b,j,1}\ne0$ implies $\b_i=0,\,j=1$. Denote
$\a'=\a+2\si_i$. Then $\a'_i=\a'_{\ol i}=0$, and (\ref{d-t-i}) with
$(m,m')=(1,0)$ becomes
\begin{equation}
D^{\a}(t^{\es_i})=\SUM{\b}{}b_{\b,1,1}t^{-(\b+\si_i)}\otimes
t^{\b+\a'-\es_{\ol i}},\mbox{ \ where }b_{\b,1,1}\ne0\mbox{ implies
}\b_i=\b_{\ol i}=0.
\end{equation}
Note that (\ref{d-t-i}) with $(m,m')=(0,1)$ becomes
\begin{equation*}
~D^{\a}(t^{\es_{\ol
i}})=\SUM{\b,j',k'}{}c^{0,1}_{\b,j',k'}t^{-(\b+j'\es_i+k'\es_{\ol
i})}\otimes t^{\b+j'\es_i+k'\es_{\ol i}+\a'-2\si_i+\es_{\ol i}} .
\end{equation*}
Applying $D^\a$ to $[t^{\es_i},t^{\es_{\ol i}}]=0$, we obtain
$t^{\es_i}\cdot D^{\a}(t^{\es_{\ol i}})=t^{\es_{\ol i}}\cdot
D^{\a}(t^{\es_i}).$ So
\begin{eqnarray}\label{iequa}
\!\!\!\!&\!\!\!\!\!\!\!\!\!\!\!\!\!\!\!\!\!\!\!\!\!\!\!\!&
\SUM{\b,j',k'}{}c^{0,1}_{\b,j',k'}\left({}^{^{^{}}}k't^{-(\b+j'\es_i+(k'+1)\es_{\ol
i})}\otimes t^{\b+j'\es_i+k'\es_{\ol i}+\a'-2\si_i+\es_{\ol
i}}\right.
\nonumber\\
\!\!\!\!&\!\!\!\!\!\!\!\!\!\!\!\!\!\!\!\!\!\!\!\!\!\!\!\!&\ \ \ \ \
\ \ \ \ \ \ \ \ \ \left.{}^{^{^{}}}+
(k'-1)t^{-(\b+j'\es_i+k'\es_{\ol i})}\otimes
t^{\b+j'\es_i+k'\es_{\ol i}+\a'-2\si_i}\right)
\nonumber\\
\!\!\!\!&\!\!\!\!\!\!\!\!\!\!\!\!\!\!\!\!\!\!\!\!\!\!\!\!&
=\SUM{\b}{}b_{\b,1,1}t^{-(\b+\si_i+\es_{\ol i})}\otimes
t^{\b+\a'-\es_{\ol i}}.
\end{eqnarray}
This forces $b_{\b,1,1}=c^{0,1}_{\b,j',k'}=0.$

\vskip5pt {\it Subcase 4.2})~~{$~\a_{\ol i}>-1$}. Applying $D^\a$ to
$[t^{\es_i},t^{m\es_{i}}]=0, m\in\Z$, we obtain
\begin{equation*}
t^{\es_i}\cdot D^\a(t^{m\es_{i}})=t^{m\es_{i}}\cdot
D^\a(t^{\es_{i}}).
\end{equation*}
We shall simplify notations by denoting
\begin{equation}\label{e-denoting}t^{\b,j,k}_{\a,\g}=t^{-(\b+j\es_i+k\es_{\ol i})}\otimes
t^{\b+j\es_i+k\es_{\ol
 i}+\a+\g}.\end{equation}
Using (\ref{d-t-i}), we obtain
\begin{eqnarray}
\label{i-2i} \nonumber
&&\SUM{\b,j,k}{}(-m(k-1)b_{\b,j+m-1,k-1}+m(k+\a_{\ol
i})b_{\b,j,k})t^{\b,j,k}_{\a,m\es_i-\es_{\ol i}}\\
&&=\SUM{\b,j,k}{}(-(k-1)c_{\b,j,k-1}^{m\es_i,0}+(k+\a_{\ol
i})c_{\b,j,k}^{m\es_i,0})t^{\b,j,k}_{\a,m\es_i-\es_{\ol i}}\,.
\end{eqnarray}
Since $D^{\a}(t^{\es_i})$ is a finite sum, we can choose $m\gg0$
such that $$b_{\b,j,k}\ne0\ \ \ \Longrightarrow\ \ \
b_{\b,j+m-1,k-1}=b_{\b,j-m+1,k+1}=0.$$ Thus (\ref{i-2i}) gives
\begin{eqnarray}
\label{m-i} &&m(k+\a_{\ol
i})b_{\b,j,k}=-(k-1)c_{\b,j,k-1}^{m\es_i,0}+(k+\a_{\ol
i})c_{\b,j,k}^{m\es_i,0},\\
&&\label{m-2i}
-mkb_{\b,j,k}=-kc_{\b,j-m+1,k}^{m\es_i,0}+(k+1+\a_{\ol
i})c_{\b,j-m+1,k+1}^{m\es_i,0}.
\end{eqnarray}
Note that $b_{\b,j,k}\neq0$ implies $k=1,-(\a_i+1)$. First assume
$b_{\b,j,1}\neq0$. Then (\ref{m-i}) gives
$$\mbox{$m(1+\a_{\ol
i})b_{\b,j,1}=(1+\a_{\ol i})c_{\b,j,1}^{m\es_i,0}$, \ \ which means
$c_{\b,j,1}^{m\es_i,0}\ne0$.} $$
 Taking $k=2,3,4,\cdots$, we get
\begin{equation}
\nonumber 0=-c_{\b,j,1}^{m\es_i,0}+(2+\a_{\ol
i})c_{\b,j,2}^{m\es_i,0}=-2c_{\b,j,2}^{m\es_i,0}+(3+\a_{\ol
i})c_{\b,j,3}^{m\es_i,0}=-3c_{\b,j,2}^{m\es_i,0}+(4+\a_{\ol
i})c_{\b,j,4}^{m\es_i,0}=\cdots.
\end{equation}
From this, we can deduce $c_{\b,j,k}^{m\es_i,0}\neq0$ for $1\le k\in
\Z$, which contradicts the fact that $D^{\a}(t^{m\es_i})$ is a
finite sum. Thus $b_{\b,j,1}=0$.

Next set $q=-(\a_{\ol i}+1)$ and assume $b_{\b,j,q}\neq0$. From
(\ref{m-2i}), we get
$$mqb_{\b,j,q}=qc_{\b,j-m+1,q}^{m\es_i,0},\mbox{ \ which implies \ }
c_{\b,j-m+1,q}^{m\es_i,0}\ne0.$$ Taking $k=q-1,q-2,\cdots$, we get
\begin{equation}
\nonumber 0=(q-1)c_{\b,j-m+1,q-1}^{m\es_i,0}+(q+\a_{\ol
i})c_{\b,j-m+1,q}^{m\es_i,0}=(q-2)c_{\b,j-m+1,q-2}^{m\es_i,0}+(q-1+\a_{\ol
i})c_{\b,j-m+1,q-1}^{m\es_i,0}=\cdots.
\end{equation}
We get a contradiction for the same reason. Thus $b_{\b,j,q}=0$.
\vskip5pt {\it Subcase 4.3})~~{$~\a_{\ol i}<-2$}. In this subcase,
since $-(\a_{\ol i}+1)>1$, we can further replace $D^\a$ by
$D^\a-u_{\rm inn}$, where $u$ has the form
$$u=\SUM{j}{}\,\SUM{k=2}{-(\a_i+1)}a_{j,k}t^{-(\b+j\es_i+k\es_{\ol i})}\otimes
t^{\b+j\es_i+k\es_{\ol i}+\a+\si_i},\mbox{ \ for some }a_{j,k}\in
\F,$$ such that $b_{\b,j,k}\neq0$ only if $k=1$ (note that such $u$
uniquely exists, and $a_{j,k}$'s are uniquely determined in order to
guarantee that after the replacement we have $b_{\b,j,k}=0$ for
$k\ne1$). Then similar to Subcase 4.2), from (\ref{m-i}), if
$mb_{\b,j,1}=c_{\b,j,1}^{m\es_i,0}\neq0$, we  get
\begin{equation}
\nonumber 0=-c_{\b,j,1}^{m\es_i,0}+(2+\a_{\ol
i})c_{\b,j,2}^{m\es_i,0}=-2c_{\b,j,2}^{m\es_i,0}+(3+\a_{\ol
i})c_{\b,j,3}^{m\es_i,0}=\cdots,
\end{equation}
and we get a contradiction again.
 \vskip5pt {\it Subcase 4.4   })~~{$~\a_{\ol
i}=-1$}. Then $b_{\b,j,k}\neq0$ only if $k=0$ or $k=1$. From
(\ref{i-i}), we  get
\begin{eqnarray}
\label{4.41}&&
(\b_i+j)b_{\b,j,0}-(\b_i+j-1)b_{\b,j-1,0}=-c^{0,1}_{\b,j,-1},\\
\label{4.42}
&&(\b_i+j)b_{\b,j,1}-(\b_i+j-1)b_{\b,j-1,1}=-c^{0,1}_{\b,j,1}.
\end{eqnarray}
If $c^{0,1}_{\b,j,-1}\neq0$ or $c^{0,1}_{\b,j,1}\neq0$, then similar
to Subcase 4.2), we can deduce $D^{\a}(t^{m\es_i})$ is an infinite
sum. Thus $c^{0,1}_{\b,j,-1}=c^{0,1}_{\b,j,1}=0$, and (\ref{4.41})
and (\ref{4.42}) become
\begin{equation}
\label{4.43} (\b_i+j)b_{\b,j,0}-(\b_i+j-1)b_{\b,j-1,0}=
(\b_i+j)b_{\b,j,1}-(\b_i+j-1)b_{\b,j-1,1}=0.
\end{equation}
Similar to Subclaim 4.1),  we deduce that if $b_{\b,j,0}\neq0$ or
$b_{\b,j,1}\neq0$, then $\b_i=0$ and $j=0$. As the result, we can
suppose
\begin{equation}\label{e-222}
d_0(t^{\es_i})=\SUM{\b}{} (b_{\b,0,0}t^{-\b}\otimes t^{\b-\es_{\ol
i}}+b_{\b,0,1}t^{-\b-\es_{\ol i}}\otimes t^\b)~.
\end{equation}
Since $\b_i=\b_{\ol i}=0$ if $\b$ appears in (\ref{e-222}), we
obtain \begin{equation} \label{e-about-i} d_0(t^{\es_i})\equiv
0\,({\rm mod\,}V^{\le i}).\end{equation}

\vskip5pt {\it Step 2.}~~We want to prove $d_0(t^{\es_{\ol
i}})\equiv0\,({\rm mod\,}V^{\le i})$.

 First, by (\ref{d-t-i}),
we can rewrite (recall notation (\ref{e-denoting}))
\begin{equation}
D^\a(t^{\es_{\ol i}})=\SUM{\b\in{\ol \G},\,j\in \Z,\,k\in
I}{}c^{0,1}_{\b,j,k}t^{\b,j,k}_{\a,\es_{\ol i}}\mbox{, \ ~where \ }
I\subset \Z \mbox{\ is a finite subset\ }.
\end{equation}
Note that (\ref{i-i}) becomes
\begin{equation}
\label{0-i-i} 0=-(\b_{\ol i}+\a_{\ol
i}+k+1)c^{0,1}_{\b,j,k}+(\b_{\ol i}+k-1)c^{0,1}_{\b,j,k-1}.
\end{equation}
Using this, discussing as in Cases 1--3) in Step 1, we obtain
\begin{equation}\label{e-----}
c^{0,1}_{\b,j,k}\ne0\ \ \ \Longrightarrow\ \ \ \b_{\ol i}=0\mbox{ \
and \ }\a_{\ol i}=\a_i\in\Z.
\end{equation}
Denote $u_{\b,j}=\sum_{k\in
I}{}c^{0,1}_{\b,j,k}t^{\b,j,k}_{\a,\si_i}$. Then (\ref{0-i-i}) means
\begin{equation}
[t^{\es_i},u_{\b,j}]=\SUM{k\in I}{}((\b_{\ol i}+\a_{\ol
i}+k+1)c^{0,1}_{\b,j,k}-(\b_{\ol
i}+k-1)c^{0,1}_{\b,j,k-1})t^{\b,j,k}_{\a,\es_i}=0.
\end{equation}
On the other hand,
\begin{equation}
\label{-inn-}[t^{\es_{\ol i}},u_{\b,j}]=\SUM{k\in I}{}((\b_
i+j)c^{0,1}_{\b,j,k}t^{\b,j+1,k}_{\a,\es_{\ol i}}-(\b_
i+j+\a_i+1)c^{0,1}_{\b,j,k}t^{\b,j,k}_{\a,\es_{\ol i}}).
\end{equation}
Thus if we replace $D^\a$ by $D^\a-u_{\rm inn}$, where $u$ is some
combination of $u_{\b,j}$, we still have (\ref{e-about-i}), and we
can further suppose
\begin{equation}
\label{e-case-+}
 ~c^{0,1}_{\b,j,k}\ne0\ \ \Longrightarrow\ \
 \left\{
\begin{array}{lll}
\mbox{1)}~~\b_{i}\ne0, &\!\!\mbox{and~~}j=0,\ \ \mbox{or}
\\[6pt]
2)~~\b_{i}=0,&\!\!\mbox{and}~~j=1 ~\mbox{or}~j=-(\a_{\ol i}+1).
\end{array}
\right. \end{equation} The rest of the proof is similar to Step 1
(but the proof is more simple than that of Step 1 since we have
(\ref{0-i-i}) instead of (\ref{i-i})).

\begin{subc}\label{subclaim--}
\addbra \rm We can further suppose $d_0(t^{2\es_i})\equiv
d_0(t^{2\es_{\ol i}})\equiv 0\,({\rm mod\,}V^{\le i})$.
\end{subc}

We shall again consider $D^\a$ for $\a\in\pi^{-1}(0)$. From the
proof of Subclaim \ref{sub2}, we can see that we only need to
consider the case when $\a_i=\a_{\ol i}\in\Z$.

 In order
to prove this subclaim, we first prove the following subclaim.
\begin{subc}\label{sub++}
\addbra \rm Suppose for some $m,{m'}\in \Z_+$,
\begin{equation}\label{condition-m-m}
t^{\es_i}\cdot(D^\a(t^{m\es_i+{m'}\es_{\ol i}}))=t^{\es_{\ol
i}}\cdot(D^\a(t^{m\es_i+{m'}\es_{\ol i}}))=0.\end{equation} Set
$q=-(\a_i+m),\,q'=-(\a_i+{m'})$. Then (\ref{d-t-i}) can be rewritten
as \begin{equation} \label{D-m-m} D^\a(t^{m\es_i+{m'}\es_{\ol
i}})=\SUM{\b}{}\ \SUM{j=q}{0}\ \SUM{k=q'}{0}\
c^{m,{m'}}_{\b,j,k}t^{\b,j,k}_{\a,m\es_i+{m'}\es_{\ol i}}\,,
\end{equation}
such that \begin{equation} \label{D-m-m+}  c^{m,{m'}}_{\b,j,k}\ne0\
\ \ \Longrightarrow\ \ \ \b_i=\b_{\ol i}=0\mbox{ \ and \ }\a_i\ge
{\rm max}\{-m,-{m'}\}.\end{equation} Furthermore,
\begin{equation}\label{e++++++}
\mbox{$D^\a(t^{m\es_i+{m'}\es_{\ol i}})=0$ \ \ if
$c^{m,{m'}}_{\b,j,k}=0$ for some $q<j<0$ or some $q'<k<0$}.
\end{equation}
\end{subc}

Using $t^{\es_i}\cdot(D^\a(t^{m\es_i+{m'}\es_{\ol i}}))=0$, for a
fixed pair $(\b,j)$, we get
\begin{equation}
\label{k-q} (\b_{\ol i}+k-1)c^{m,{m'}}_{\b,j,k-1}=(\b_{\ol
i}+k+\a_{\ol i}+{m'})c^{m,{m'}}_{\b,j,k}\, \mbox{\ for  \ } k\in \Z.
\end{equation}
First suppose $\b_{\ol i}\neq0$. Then $\b_{\ol i}+k-1\neq0$ for
$k\in \Z$ by (\ref{e-require}). Thus (\ref{k-q}) shows that if
$c^{m,{m'}}_{\b,j,k}\neq0$ for one $k\in\Z$ then
$c^{m,{m'}}_{\b,j,k}\ne0$ for all $k\in\Z$. Thus
$c^{m,{m'}}_{\b,j,k}\ne0$ implies $\b_{\ol i}=0$.

Next suppose $c^{m,{m'}}_{\b,j,k}\neq0$ for some $k>0$ or
$k<-(\a_{\ol i}+{m'})$ or $\a_{\ol i}<-{m'}$. Then (\ref{k-q}) can
again deduce $c^{m,{m'}}_{\b,j,k}\neq0$ for infinite many $k$. Thus
$c^{m,{m'}}_{\b,j,k}\neq0$ implies $-(\a_{\ol i}+{m'})\le k\le0$ and
$\a_{\ol i}\ge-{m'}$. Furthermore, if $c^{m,{m'}}_{\b,j,k}=0$ for
some $q'<k<0$, one can deduce from (\ref{k-q}) that
$c^{m,{m'}}_{\b,j,k}=0$ for all $q'<k<0$. Using these results and
the symmetric results derived from $t^{\es_{\ol
i}}\cdot(D^\a(t^{m\es_i+{m'}\es_{\ol i}}))=0$, we obtain
(\ref{D-m-m})--(\ref{e++++++}). \vskip5pt
 Now we prove Subclaim \ref{subclaim--}).
First by Subclaim \ref{sub2}), we see that condition
(\ref{condition-m-m}) holds for $(m,m')=(2,0),\,(0,2)$. Thus by
Subclaim \ref{sub++}), we can suppose
\begin{equation}
\nonumber D^\a(t^{2\es_i})=\SUM{\b}{}\ \SUM{j=-(\a_i+2)}{0}\
\SUM{k=-\a_i}{0}\ c^{2,0}_{\b,j,k}t^{\b,j,k}_{\a,2\es_i},~~~~\ \ \
D^\a(t^{2\es_{\ol i}})=\SUM{\b}{}\ \SUM{j=-\a_i}{0}\
\SUM{k=-(\a_i+2)}{0}\ c^{0,2}_{\b,j,k}t^{\b,j,k}_{\a,2\es_{\ol i}}~.
\end{equation}
Applying $D^\a$ to $[t^{2\es_i},t^{2\es_{\ol i}}]=4t^{\si_i}$, by
Claim \ref{clai-sub1}, we get, for a fixed $\b$,
\begin{eqnarray}
\nonumber && \SUM{j=-(\a_i+2)}{0}\ \SUM{k=-\a_i}{0}\
(-2(j-1)c^{2,0}_{\b,j-1,k+1}+(j+\a_i+2)c^{2,0}_{\b,j,k})
t^{\b,j,k}_{\a,2\si_i}\\
\nonumber &&=\SUM{j=-\a_i}{0}\ \SUM{k=-(\a_i+2)}{0}\
(-2(k-1)c^{0,2}_{\b,j+1,k-1}+(k+\a_{\ol i}+2)c^{0,2}_{\b,j,k})
t^{\b,j,k}_{\a,2\si_i}.
\end{eqnarray}
Taking $j=-(\a_i+1),k=-\a_i$, we get
$c^{2,0}_{\b,-(\a_i+1),-\a_i}=0$. Taking $j=-\a_i,k=-(\a_i+1)$, we
get $c^{0,2}_{\b,-\a_i,-(\a_i+1)}=0$. By (\ref{e++++++}), we obtain
Subclaim
 \ref{subclaim--}).

\begin{subc}\addbra
\rm We can further suppose $d_0(t^{\si_i+\es_i})\equiv
d_0(t^{\si_i+\es_{\ol i}})\equiv 0\,({\rm mod\,}V^{\le i})$.
\end{subc}

Again we consider $D^\a$ with $\a_i=\a_{\ol i}\in\Z$. Note that
$t^{\si_i+\es_i}$ and $t^{\si_i+\es_{\ol i}}$ correspond
respectively to $(m,m')=(2,1)$ and $(1,2)$ in Subclaim \ref{sub++}),
and condition (\ref{condition-m-m}) holds for these two cases. Thus
we can suppose $\a_i\ge2$.

 Applying $D^\a$ to
$[t^{2\es_i},[t^{2\es_{\ol i}},t^{\si_i+\es_i}]=-8t^{\si_i+\es_i}$,
we get
\begin{equation}
\nonumber t^{2\es_i}\cdot t^{2\es_{\ol i}}\cdot
D^\a(t^{\si_i+\es_i})=-8D^\a(t^{\si_i+\es_i}).
\end{equation}
Using (\ref{D-m-m}),  we obtain
\begin{equation}
\nonumber
(k-1)(jc^{2,1}_{\b,j,k}-(j+\a_i+3)c^{2,1}_{\b,j+1,k-1})-(k+\a_i+2)((j-1)
c^{2,1}_{\b,j-1,k+1}-(j+\a_i+2)c^{2,1}_{\b,j,k})=2c^{2,1}_{\b,j,k}.
\end{equation}
Noting from (\ref{D-m-m}) that $c^{2,1}_{\b,j,k}=0$ for $j>0$ or
$k>0$. Taking $j=0,$ we get
$$(k+\a_i+2)c^{2,1}_{\b,-1,k+1}=\a_ic^{2,1}_{\b,j,k}.$$ Taking
$k=0$, we deduce
$c^{2,1}_{\b,-1,0}=\a_i^{-1}(\a_i+2)c^{2,1}_{\b,-1,1}=0$. Thus by
(\ref{e++++++}), we get $D^\a(t^{\si_i+\es_{i}})=0$ . Similarly,
$D^\a(t^{\si_i+\es_{\ol i}})=0$.
\begin{subc}\addbra
\rm We can suppose that $d_0(t^{m\es_i+{m'}\es_{\ol i}})\equiv
0\,({\rm mod\,}V^{\le i})$ for all $m,{m'}\in\Z$.
\end{subc} \ \indent
Denote $$\PP_i^+={\rm span}\{t^{j\es_i+k\es_{\ol i}}\,|\,1\le
j,k\in\Z\},$$ a subalgebra of $\PP$. It is easy to verify that $A=\{
t^{\es_p},t^{2\es_p},t^{\si_p+\es_p}|p=i,{\ol i} \}$ is a generating
set of $\PP_i^+$. So $d_0(x)\equiv 0\,({\rm mod\,}V^{\le i})$ for
$x\in\PP_i^+$. For any $m,{m'}\in\Z$, we can choose $\wt m,{\wt
m'}\in\Z^+$ such that $m+\wt m\gg0,\,{m'}+{\wt m'}\gg0$. Applying
$d_0$ to the equation
\begin{equation}
\nonumber [t^{m\es_i+{m'}\es_{\ol i}},t^{\wt m\es_i+{\wt m'}\es_{\ol
i}}]=(m{\wt m'}-{m'}\wt m)t^{(m+m')\es_i+({\wt m}+{\wt m'})\es_{\ol
i}{\ssc\,}-\si_i}\in\pp_i^+,
\end{equation}
we get $t^{\wt m\es_i+{\wt m'}\es_{\ol i}}\cdot
d_0(t^{m\es_i+{m'}\es_{\ol i}})=0$ for all $\wt m,{\wt m'}\gg0$.
Similar to the proof of Lemma \ref{lemma2.3}, we get
$d_0(t^{m\es_i+{m'}\es_{\ol i}})\equiv 0\,({\rm mod\,}V^{\le i})$.
\vskip5pt Now we shall finish the proof of (\ref{e-ind}) (namely,
Claim \ref{claim3}). Let $\g\in\G$. If $\g_i=\g_{\ol i}=0$, then by
applying $d_0$ to the equations
\begin{equation}
\nonumber [t^{m\es_i},t^\g]=0,\ \ ~ [t^{m\es_{\ol i}},t^\g]=0\mbox{
\ \ for \ }m\in\Z,
\end{equation}
we get $t^{m\es_i}\cdot d_0(t^\g)\equiv t^{m\es_{\ol i}}\cdot
d_0(t^\g)\equiv 0\,({\rm mod\,}V^{\le i})$. By Lemma \ref{lemma2.3},
we get $d_0(t^\g)\equiv 0\,({\rm mod\,}V^{\le i})$.

Suppose  $(\g_i,\g_{\ol i})\neq0$. Simply write $d_0(t^\g)$ as
\begin{equation}\label{d-0-g}d_0(t^\g)=\SUM{\a,\b}{}\,c^\g_{\a,\b}t^\a\otimes
t^\b\mbox{ \ \ for some \ }c^\g_{\a,\b}\in\F.\end{equation} For
convenience, we denote
$$a_{j,k}=-((j+2)\g_{\ol i}-(k+2)\g_i)(j(k+1+\g_{\ol
i})-k(j+1+\g_i)).$$ Applying $d_0$ to
\begin{equation}
\nonumber [t^{-j\es_i-k\es_{\ol i}},[t^{(j+2)\es_i+(k+2)\es_{\ol
i}},t^\g]]=a_{j,k}t^\g,
\end{equation}
we obtain
\begin{equation}
\nonumber t^{-j\es_i-k\es_{\ol i}}\cdot t^{(j+2)\es_i+(k+2)\es_{\ol
i}}\cdot \SUM{\a,\b}{}\,c^\g_{\a,\b}t^{\a}\otimes
t^\b=a_{j,k}\SUM{\a,\b}{}\,c^\g_{\a,\b}t^{\a}\otimes t^\b.
\end{equation}
This means
\begin{eqnarray*}
\!\!\!\!\!\!\!\!&\!\!\!\!\!\!\!\!\!\!\!\!\!\!\!\!&
a_{j,k}\SUM{\a,\b}{}\,c^\g_{\a,\b}t^{\a}\otimes t^\b\\
\!\!\!\!\!\!\!\!&\!\!\!\!\!\!\!\!\!\!\!\!\!\!\!\!& =
-\SUM{\a,\b}{}\,(c^\g_{\a,\b}((j+2)\a_{\ol
i}-(k+2)\a_i)(j(\a_{\ol i}+k+1)-k(\a_i+j+1))t^\a\otimes t^\b\\
\!\!\!\!\!\!\!\!&\!\!\!\!\!\!\!\!\!\!\!\!\!\!\!\!&
+c^\g_{\a,\b}((j+2)\a_{\ol i}-(k+2)\a_i)(j\b_{\ol i}-k\b_i)
t^{\a+(j+1)\es_i+(k+1)\es_{\ol i}}\otimes t^{\b-(j+1)\es_i-(k+1)\es_{\ol i}}\\
\!\!\!\!\!\!\!\!&\!\!\!\!\!\!\!\!\!\!\!\!\!\!\!\!&
+c^\g_{\a,\b}((j+2)\b_{\ol i}-(k+2)\b_i)(j\a_{\ol i}-k\a_i)
t^{\a-(j+1)\es_i-(k+1)\es_{\ol i}}\otimes t^{\b+(j+1)\es_i+(k+1)\es_{\ol i}}\\
\!\!\!\!\!\!\!\!&\!\!\!\!\!\!\!\!\!\!\!\!\!\!\!\!&
+c^\g_{\a,\b}((j+2)\b_{\ol i}-(k+2)\b_i)(j(\b_{\ol
i}+k+1)-k(\b_i+j+1)))t^\a\otimes t^\b.
\end{eqnarray*}
Choosing $j,k\gg0$, and comparing coefficients in the both sides, we
obtain
\begin{eqnarray*}
c^\g_{\a,\b}((j+2)\a_{\ol i}-(k+2)\a_i)(j\b_{\ol i}-k\b_i)
=c^\g_{\a,\b}((j+2)\b_{\ol i}-(k+2)\b_i)(j\a_{\ol i}-k\a_i)=0.
\end{eqnarray*}
This shows that \begin{equation}\label{c-a-b}c^\g_{\a,\b}\ne0\ \ \
\Rar\ \ \ \a_i=\a_{\ol i}=0\mbox{ \ or \ }\b_i=\b_{\ol
i}=0.\end{equation} Now we need to consider the following two cases:
$$\mbox{(i) \ }
(\a_i,\a_{\ol i})\ne0\mbox{ and }t^\b\notin \pp^{\le i};\ \ \ \ \
\mbox{(ii) \ }(\b_i,\b_{\ol i})\ne0\mbox{ and }t^\a\notin \pp^{\le
i}.$$ Without loss of generality, we can suppose for some $\g\in\G$,
there exists one term $t^\xi\otimes t^\eta$ such that
$\xi_{\ol{i}}\neq0,\xi_{\ol{j}}=\xi_j=0$ for some $j<i$ and
$\eta_i=\eta_{\ol{i}}=0$. Write
\begin{equation}
d_0(t^\g)=ct^\xi\otimes t^\eta
+\SUM{\a,\b}{}\,c^\g_{\a,\b}t^\a\otimes t^\b,
\end{equation}
where $0\neq c\in \F$. From assumption we know  there exists
$\eta_j\neq0$ or $\eta_{\ol{j}}\neq0$ for some $j<i$. Fixing $j$
and set
$d_0(t^{\es_{\ol{i}}+\es_j})=\sum_{\tau,\tau'}c_{\tau,\tau'}t^\tau\otimes
t^{\tau'}$, then from
$[t^{k\es_i},t^{\es_{\ol{i}}+\es_j}]=kt^{(k-1)\es_i+\es_j}$ we get
\begin{equation}
k\,d_0(t^{(k-1)\es_i+\es_j})=t^{k\es_i}\cdot
d_0(t^{\es_{\ol{i}}+\es_j}).
\end{equation}
For any $k\neq0$, applying $d_0$ to
$[t^{(k-1)\es_i+\es_j},t^\g]=\g_{\ol{j}}t^{\g-\es_{\ol{j}}}
+(k-1)\g_{\ol{i}}t^{\g+(k-1)\es_i-\si_{i}}$ we have
\begin{equation}
\label{ffff} t^{(k-1)\es_i+\es_j}\cdot d_0(t^\g)-t^\g\cdot
d_0(t^{(k-1)\es_i+\es_j})=d_0(\g_{\ol{j}}t^{\g-\es_{\ol{j}}}
+(k-1)\g_{\ol{i}}t^{\g+(k-1)\es_i-\si_{i}}).
\end{equation}
Since
\begin{eqnarray*}
&&t^{(k-1)\es_i+\es_j}\cdot d_0(t^\xi\otimes
t^\eta)=(k-1)(\xi_{\ol{i}}t^{(k-1)\es_i+\es_j+\xi-\si_i}\otimes
t^\eta+\eta_{\ol{i}}t^\xi\otimes
t^{(k-1)\es_i+\es_j+\eta-\si_i})~~~~~\\
&&t^\g\cdot
d_0(t^{(k-1)\es_i+\es_j})=\SUM{\tau,\tau'}{}c_{\tau,\tau'}
k(\tau_{\ol{i}}(\g_i(\tau_{\ol{i}}-1)-\g_{\ol{i}}(\tau_i+k-1))
t^{\tau+k\es_i+\g-2\si_i}\otimes t^{\tau'}\\
&&~~~~~~~~~~~~~~~~~~~~~~~~~+\tau_{\ol{i}}'(\g_i(\tau_{\ol{i}}'-1)-
\g_{\ol{i}}(\tau_i'+k-1)))t^\tau\otimes
t^{\tau'+k\es_i+\g-2\si_i}+\SUM{\mu,\nu}{}c_{\mu,\nu}t^\mu\otimes
t^\nu,
\end{eqnarray*}
where $\sum_{\mu,\nu}c_{\mu,\nu}t^\mu\otimes t^\nu\in V,$ we know
the coefficients of the term $t^{(k-1)\es_i+\es_j+\xi-\si_i}\otimes
t^\eta$ is not equal to zero, but from the discussion above, we know
$t^{(k-1)\es_i+\es_j+\xi-\si_i}\otimes t^\eta$ does not appear to
the right-hand side of equation (\ref{ffff}). This is a
contradiction. Thus we know $d_0(t^\g)\equiv 0\,({\rm mod\,}V^{\le
i})$. This proves (\ref{e-ind}) and Claim \ref{claim3}.

\begin{clai}\adddot \rm $d_\a=0$ for
all but a finite number of $\a\in\G$.
\end{clai}
\ \indent Let $u_\a\in V_\a$ such that $d_\a=(u_\a)_{\rm inn}$ for
$\a\in\G$. If $u_\a\ne0$ for infinite many $\a\in\G$, then there
exists some $1\le i\le n$ such that $(\a_i,\a_{\ol i})\ne(0,0)$ for
infinite many $\a$ with $u_\a\ne0$. So we can shoose
$t^{j\es_i+k{\es_{\ol i}}}$ such that $j\a_{\ol i}-k\a_i\ne0$ for
infinite many $\a$  with $u_\a\ne0$. Then $t^{j\es_i+k{\es_{\ol
i}}}\cdot u_\a\ne0$ (cf.~(\ref{braket})) for infinite many $\a$ with
$u_\a\ne0$, and so $d(t^{j\es_i+k{\es_{\ol
i}}})=\sum_{\a\in\G}t^{j\es_i+k{\es_{\ol i}}}\cdot u_\a$ is an
infinite sum, a contradiction. This proves the claim and Theorem
\ref{lemma3?}. \QED

\vskip6pt \ni{\it Proof of Theorem \ref{main}(1).} Let
$(\PP,[\cdot,\cdot],\D)$ be a Lie bialgebra structure on $\PP$. By
(\ref{bLie-d}), (\ref{deriv}) and Theorem \ref{lemma3?}, $\D=\D_r$
is defined by (\ref{D-r}) for some $r\in\PP\otimes\PP$. By
(\ref{cLie-s-s}), ${\rm Im}\,\D\subset{\rm Im}(1\otimes1-\tau)$.
Thus by Lemma \ref{last}, $r\in{\rm Im}(1\otimes1-\tau)$.
Then  (\ref{cLie-j-i}), (\ref{add-c}) and Theorem \ref{main}(2) show
that $c(r)=0$. Thus Definition \ref{def2} says that
$(\PP,[\cdot,\cdot],\D)$ is a triangular Lie bialgebra. \QED
\vs{10pt}

\vskip12pt

\cl{\bf References}\vs{0pt}

\vskip5pt\small
\parindent=8ex\parskip=2pt\baselineskip=2pt
\re{D1} V.G. Drinfeld, Hamiltonian structures on Lie groups, Lie
bialgebras and the geometric meaning of classical Yang-Baxter
equations, (Russian) {\it Dokl. Akad. Nauk SSSR} {\bf268} (1983),
285--287.

\re{D2} V.G. Drinfeld, Quantum groups, in: {\it Proceeding of the
International Congress of Mathematicians}, Vol.~1, 2, Berkeley,
Calif.~1986, Amer.~Math.~Soc., Providence, RI, 1987, pp.~798--820.

\re{Dzh} A.S. Dzhumadil'daev, Quasi-Lie bialgebra structures of
$sl_2$, Witt and Virasoro algebras, in : {\it Quantum deformations
of algebras and their representations} (Ramat-Gan, 1991/1992;
Rehovot, 1991/1992), 13-24, Israel Math. Conf. Proc., {\bf7},
Bar-Ilan Univ., Ramat Gan, 1993.

 \re{DZ} D. Dokovic, K. Zhao, Derivations, isomorphisms and
second cohomology of generalized Witt algebras, {\it Trans. Amer.
Math. Soc.} {\bf350} (1998), 643--664.


\re{M} W. Michaelis, A class of infinite-dimensional Lie bialgebras
containing the Virasoro algebra, {\it Adv. Math.} {\bf107} (1994),
365--392.

%

\re{NT} S.-H. Ng, E.J. Taft, Classification of the Lie
bialgebra structures on the Witt and Virasoro algebras,
{\it J. Pure Appl.~Algebra} {\bf151} (2000), 67--88.

\re{N} W.D. Nichols, The structure of the dual Lie coalgebra
of the Witt algebra, {\it J. Pure Appl. Algebra} {\bf 68} (1990),
395--364.


\re{OZ}  J.M. Osborn, K. Zhao, Generalized Poisson brackets and Lie
algebras for type $H$ in characteristic 0, {\it Math. Z.} {\bf230}
(1999), 107--143.

\re{S} G.~Song, The structure of infinite dimensional non-graded Lie
algebras and Lie superalgebras of $W$-type and related problems,
{\it Ph.~D.~Thesis,} Shanghai Jiaotong University (2005).

\re{SS} G. Song, Y. Su, Lie bialgebras of generalized Witt type,
{\it Science in China A} {\bf49} (2006), 533--544.

\re{Su} Y. Su, Poisson brackets and structure of nongraded
Hamiltonian Lie algebras related to locally-finite derivations, {\it
Canad. J. Math.} {\bf 55} (2003), 856--896.

\re{SX} Y. Su, X. Xu, Central simple Poisson algebras, {\it Science
in China A} {\bf 47} (2004), 245--263.



\re{T} E.J. Taft, Witt and Virasoro algebras as Lie
bialgebras, {\it J. Pure Appl.~Algebra} {\bf87} (1993), 301--312.

\re{X} X. Xu, New generalized simple Lie algebras of Cartan
type over a field with characteristic 0,  {\it J.~Algebra} {\bf 224}
(2000), 23--58.

\re{WSS} Y. Wu, G. Song, Y. Su,
 Lie bialgebras of generalized Virasoro-like type, {\it Acta Mathematica Sinica, English Series},
in press.

\end{document}